\documentclass[12pt]{article}

\usepackage{amssymb,amsthm,amsmath}
\usepackage{color}
\usepackage{multirow}
\usepackage{xcolor}
\usepackage{eucal}
\usepackage[linesnumbered,lined,boxed,commentsnumbered]{algorithm2e}
\usepackage{tikz}

\newcommand{\comment}[1]{}

\newcommand{\tw}{{\mathrm{tw}}}
\newcommand{\GP}{{\mathrm{GP}}}

\newcommand{\NP}{$\mathsf{NP}$}
\newcommand{\EXPTIME}{$\mathsf{EXPTIME}$}

\newtheorem{theorem}{Theorem}[section]
\newtheorem{lemma}[theorem]{Lemma}
\newtheorem{corollary}[theorem]{Corollary}

\theoremstyle{definition}

\newcommand{\bigO}{\CMcal{O}}

\def \cC {{\cal C}}

\renewcommand{\geq}{\geqslant}
\renewcommand{\leq}{\leqslant}
\tikzstyle{vertex}=[auto=left,circle,fill=none,draw=black,,minimum size=11pt, inner sep=0pt]

\title{Cops that Surround a Robber}

\author{Andrea C.~Burgess\thanks{Department of Mathematics and Statistics, University of New Brunswick, Saint John, NB, Canada. andrea.burgess@unb.ca},
Rosalind A.~Cameron\thanks{School of Mathematics and Statistics, University of Canterbury, Christchurch, New Zealand. rosalind.cameron@canterbury.ac.nz},
Nancy E.~Clarke\thanks{Department of Mathematics and Statistics, Acadia University, Wolfville, NS, Canada. nancy.clarke@acadiau.ca},
\\
Peter Danziger\thanks{Department of Mathematics, Ryerson University, Toronto, ON, Canada. danziger@ryerson.ca},
Stephen Finbow\thanks{Department of Mathematics and Statistics, St.\ Francis Xavier University, Antigonish, NS, Canada. sfinbow@stfx.ca},
Caleb W.~Jones\thanks{Department of Mathematics and Statistics, Memorial University of Newfoundland, St.~John's, NL, Canada. cwjones@mun.ca},
David A.~Pike\thanks{Department of Mathematics and Statistics, Memorial University of Newfoundland, St.~John's, NL, Canada. dapike@mun.ca}
}

\begin{document}
\maketitle

\begin{abstract}
We introduce the game of Surrounding Cops and Robbers on a graph, as a variant of the original
game of Cops and Robbers.  In contrast to the original game in which the cops win by occupying the same vertex as the robber,
they now win by occupying each of the robber's neighbouring vertices.
We denote by $\sigma(G)$ the {\em surrounding cop number} of $G$, namely the least number of cops required to surround a robber in the graph $G$.
We present a number of results regarding this parameter,
including general bounds as well as exact values for several classes of graphs.
Particular classes of interest include product graphs, graphs arising from combinatorial designs, and generalised Petersen graphs.
\end{abstract}

\section{Introduction}

The vertex-pursuit game of Cops and Robbers was introduced about 40 years ago by
Nowakowski and Winkler~\cite{NW1983} and by Quilliot~\cite{Quilliot}.
Initially only a single cop and a single robber were considered, but
in 1984 Aigner and Fromme generalised the game to include a team of cops as one of the two players in the game~\cite{AignerFromme1984};
it is this version of the game that is now the standard scenario.
Under its rules the game commences with the cops choosing their positions at vertices of a graph $G$,
next the robber chooses a starting position at a vertex, and thereafter the cops and robber alternate their respective moves.
During a move each individual is permitted to remain stationary or to take up a new vertex position that is adjacent to the individual's current position.
If a cop ever occupies the same vertex as the robber, the cops win the game.  Otherwise, if the robber can forever evade the cops, the robber is deemed to win.
The {\em cop number} $c(G)$ denotes the least number of cops for which the cops always have a winning strategy for the graph $G$.
A substantial body of literature has been written about the game of Cops and Robbers, including a number of variants of the game.
For a detailed survey we recommend the monograph by Bonato and Nowakowski~\cite{BonatoNowakowski}.  Following~\cite{BonatoNowakowski}, we will use the pronouns she/her to refer to an individual cop and he/him/his to refer to the robber.

In this paper, we introduce a new variation of the game. The players are a set of $k>0$ cops and a single robber.
Again, the cops choose their starting vertex positions, the robber then chooses an initial position on a vertex unoccupied by a cop,
and the cops and robber thereafter alternate moves.  On each move, a player may move to an adjacent vertex or pass, subject to the following restrictions:  the robber is forbidden to move to or stay at a vertex that is occupied by a cop.
If, at any time, the robber finds that his position is also occupied by a cop, then he is immediately compelled to move to a neighbouring vertex not occupied by a cop.
Having a cop occupy the same vertex as the robber is no longer how the game is won by the cops.
Instead, the cops win the game if at any time
each of the neighbours of the robber's vertex is occupied by a cop;
this would also be when the cops proclaim {\em ``Surrender!  You're Surrounded!''}
If the robber can forever avoid the circumstance of being surrounded, then the robber wins.
The game is played with perfect information (that is, each player knows the location of every player)
and the cops are able to communicate and coordinate their moves.
By $\sigma(G)$ we denote the least number of cops for which the cops always have a winning strategy.
We call $\sigma(G)$ the {\em surrounding cop number} for the graph $G$.

Observe that any strategy by which the cops can successfully surround the robber
yields a strategy by which the cops can also win the original game.  Hence the surrounding cop number is an upper bound on the cop number.

\begin{lemma}
\label{Lemma-copVsurround}
For any graph $G$, $c(G) \leq \sigma(G)$.
\end{lemma}

By $K_{1,n-1}$, $P_n$ and $C_n$ we denote the star, path and cycle, respectively, on $n$ vertices.
For stars and most cycles, the bound of Lemma~\ref{Lemma-copVsurround} is an equality
(the only exception is the 3-cycle, for which $c(C_3)=1$ and $\sigma(C_3)=2$).
For paths, observe that $\sigma(P_n) = 1$ if $n \leq 3$ and $\sigma(P_n) = 2$ when $n \geq 4$, whereas $c(P_n)=1$ for all $n$.
More generally, if $G$ is a tree other than a star or a path then $c(G)=1$ and $\sigma(G)=2$.
If $G$ is a wheel graph (the graph obtained by joining a vertex $x$ to each vertex of a cycle $C_{n-1}$), then again $c(G)=1$ but now $\sigma(G) = 3$.

In order for the cops to surround the robber at some vertex, clearly the number of cops must be at least the degree of the robber's vertex.
Hence the minimum degree $\delta(G)$ of the graph $G$ is a lower bound on the surrounding cop number.

\begin{lemma}
\label{Lemma-deltabound}
For any graph $G$, $\delta(G) \leq \sigma(G)$.
\end{lemma}

Another straightforward bound is based on the clique number of $G$, $\omega(G)$ (the maximum number of pairwise adjacent vertices in $G$).

\begin{lemma}\label{Lemma-OmegaMinusOne}
For any graph $G$, $\omega(G)-1 \leq \sigma(G)$.
\end{lemma}

\begin{proof}
If there are at most $\omega(G)-2$ cops, the robber initially selects a vertex within a maximum clique.
As there are too few cops to surround the robber on this vertex, he moves from it only if a cop moves onto his position.
When this occurs, however,  there will be at least one vacant vertex of the same clique that he can then occupy.
\end{proof}

We note that these three lower bounds on $\sigma(G)$ are tight.  For instance, $\sigma(K_n) = \delta(K_n) = n-1 = \omega(K_n)-1$ and, for $n \geq 4$, $\sigma(C_n) = 2 = c(C_n)$.

Although in cases such as when $G$ is an $n$-cycle and $n \geq 4$ the two parameters $c(G)$ and $\sigma(G)$ are equal in value,
the two parameters are, in general, very distinct from one another.
The difference $\sigma(G) - c(G)$ can be arbitrarily large, as illustrated by complete graphs, for which $c(K_n) = 1$ and $\sigma(K_n) = n-1$,
and complete bipartite graphs $K_{m,n}$ with $2 \leq m \leq n$, for which $c(K_{m,n})=2$ and $\sigma(K_{m,n})=m$.

With respect to the bound of Lemma~\ref{Lemma-deltabound},
the difference between $\sigma(G)$ and $\delta(G)$ can also be arbitrarily large, such as when
$G$ consists of a large clique with one additional vertex attached via a single pendant edge.
Regarding Lemma~\ref{Lemma-OmegaMinusOne},
we will see later (in Theorem~\ref{Thm-StrongProductStar}) an example of a family of graphs which shows that the difference between $\sigma(G)$
and the lower bound of $\omega(G)-1$ can also be arbitrarily large.

An elementary upper bound on $\sigma(G)$ involves the independence number $\alpha(G)$
(the maximum number of pairwise non-adjacent vertices in $G$).

\begin{lemma}\label{Lemma-UpperBound}
For any graph $G$, $\sigma(G) \leq |V(G)| - \alpha(G)$.
\end{lemma}

\begin{proof}
Let $S$ be any vertex cover in $G$, the complement of which is an independent set.
If $v \in V(G) \setminus S$ then $N(v) \subseteq S$ and so initially placing a cop at each vertex of $S$ is a trivial winning strategy for the cops.
The cardinality of a minimum vertex cover is $|V(G)| - \alpha(G)$.
\end{proof}

Combining Lemma~\ref{Lemma-UpperBound} with Lemma~\ref{Lemma-OmegaMinusOne} yields the following result,
which was previously observed by Chartrand and Schuster in the context of \linebreak Nordhaus-Gaddum relations~\cite{CS1974}.

\begin{corollary}
For any graph $G$, $\alpha(G) + \omega(G) \leq |V(G)| + 1$.
\end{corollary}

Before proceeding to delve further into the surrounding cop number, we pause here
to briefly comment on two Cops and Robbers variants which at first glance appear similar to Surrounding Cops and Robbers. In the game of Containment~\cite{CKM2020,Pralat2015}, the robber plays on vertices while the cops play on edges; the cops win if they occupy all of the robber's incident edges. In Containment, however, the cops are unable to force the robber to move from a vertex; thus the maximum degree is a lower bound on the containability number, whereas the minimum degree is a lower bound for $\sigma(G)$.
For several graphs, including the Petersen Graph and the Cartesian product of paths, the containability number and surrounding cop number differ, further illustrating the distinction between these variants.

In the game of Cheating Robot Cops and Robbers~\cite{HugganNowakowski}, the cops and robber play on vertices and move simultaneously on a round; however the robber is a robot which cheats and hence knows in advance what the cops' move will be. The cops capture the robot if both players either traverse the same edge in a round, or land on the same
vertex at the end of the round. Any winning strategy for the cops requires them to surround the robot's position.
Viewed as a game with sequential turns, the robot is compelled to move from a vertex immediately after a cop lands on it, as
in the Surrounding Cops and Robbers game. However, since traversing the same edge as a cop in a given round would cause the robot to lose,
the robot is thus forbidden from moving to the vertex which the cop occupied in the previous
round, a move that is allowed in Surrounding Cops and Robber.  This added restriction is significant; for example, on a path with at least four vertices, $\sigma(P_n) = 2$ but only one cop is required to capture a cheating robot.

To provide a brief outline of the remainder of this paper, Section~\ref{Section-Preliminaries} includes additional bounds on the parameter as well as other pertinent observations.  An algorithm that determines the surrounding cop number of a given graph is presented in Section~\ref{Section-Algorithmic}.  In Section~\ref{Section-GraphClasses}, results about and bounds on the surrounding cop number are established for
a number of graph classes, such as graphs arising from product operations, graphs based on combinatorial designs, and generalised Petersen graphs.  We conclude with some further discussion, questions, and open problems.

\section{Additional bounds and other observations}
\label{Section-Preliminaries}

In this section we present some less obvious bounds on the surrounding cop number, as well as some other useful observations about this parameter.

As defined in~\cite{BonatoNowakowski}, $H$ is a {\em retract} of $G$ if there is a homomorphism $f:V(G)\rightarrow V(H)$ so that $f(x)=x$ for all $x\in V(H)$.
Retracts supply a lower bound on $c(G)$ {\cite{Berarducci93}}.  A similar proof provides an analogous lower bound on the surrounding cop number of a graph.

\begin{theorem}\label{Theorem-retract}
If $H$ is a retract of $G$ then $\sigma(H)\leq \sigma(G)$.
\end{theorem}

In the literature, lower bounds have been given on the cop number of graphs of sufficiently large girth.  In~\cite{AignerFromme1984}, it is proved that any graph $G$ of girth at least 5 satisfies $c(G) \geq \delta(G)$.  Frankl~\cite{Frankl1987} extended this result to prove that if the girth of $G$ is at least $8t-3$ (where $t$ is a positive integer) and $\delta(G)>d$, then $c(G)>d^t$.
In~\cite{CKM2020} it is proved that for all $\delta \geq 3$, every $\delta$-regular graph $G$ with girth at least 5 has containability number $\xi(G) \geq \delta(G)+1$ and every $\delta$-regular graph $G$ with girth at least 7 has containability number $\xi(G) \geq \delta(G)+2$.  Here we give a lower bound on the surrounding cop number of graphs with girth at least 7.

\begin{theorem} \label{Theorem-Girth}
If $G$ has girth $g \geq 7$ and minimum degree $\delta \geq 3$, then $\sigma(G) \geq \delta+1$.
\end{theorem}

\begin{proof}
Suppose there are $\delta$ cops.

We first show that after the cops have chosen their initial vertices, the robber may choose a vertex which guarantees he will not be surrounded after the cops' next move.  Assuming that this is not possible, it follows that every cop must be within distance at most two from every vertex not containing a cop.  It is clear that there must be a vertex $R$ which is not adjacent to every cop, as otherwise the girth of $G$ is 3 or 4. The robber begins the game on vertex $R$, and is not surrounded prior to the cops' next move.  Let $C$ be a vertex containing a cop which is not adjacent to $R$, and let $x$ be a common neighbour of $C$ and $R$.   Since $G$ has girth greater than 4, $x$ is the only common neighbour of $C$ and $R$.  Since the cops will surround $R$ on the next turn, it follows that there must be a cop on a vertex $C' \neq C$.  Moreover, $x$ cannot contain a cop, as any cops on $C$ and $x$ cannot reach distinct neighbours of $R$ on their next turn.  Thus  $C'$ must be adjacent to either $x$ or a neighbour of $x$, as well as to either $R$ or a neighbour of $R$.  But as $x$ and $R$ are adjacent, $G$ must contain a cycle of length at most 5, contradicting that its girth is at least 7.

It remains to show that following the initial round, after any move by the cops in which the robber is not surrounded, the robber has a move which will guarantee that he is not surrounded after the cops' next move.  Let $R$ be the robber's vertex, and let $X=\{x_1, \ldots, x_{\deg(R)}\}$ be the set of neighbours of $R$.

For $i \in \{1, 2, \ldots, \deg(R)\}$, let $S_i$ denote the set of $\deg(x_i)-1$ vertices other than $R$ which are adjacent to $x_i$, and let $T_i$ denote the set of vertices other than $x_i$ which are adjacent to a vertex in $S_i$.  Note that, since the girth of $G$ is at least 7, we have that for any $i,j \in \{1, \ldots, \deg(R)\}$ with $i \neq j$:
\begin{enumerate}
\item $x_j \notin S_i \cup T_i$ (otherwise $G$ would contain a $3$- or $4$-cycle)
\item $S_i \cap S_j = \emptyset$ (otherwise $G$ would contain a $4$-cycle)
\item $S_i \cap T_j = \emptyset$ (otherwise $G$ would contain a $5$-cycle)
\item $T_i \cap T_j = \emptyset$ (otherwise $G$ would contain a $6$-cycle).
\end{enumerate}
In particular, $S_i \cup T_i$ is disjoint from $\{x_j\} \cup S_j \cup T_j$ whenever $i \neq j$.

First suppose that one of the neighbours, say $x_i$, of $R$ is such that $N[x_i]$ contains no cop.  In this case, the robber passes.  (Note that, in particular, $R$ contains no cop and so the robber is not forced to move.)  On the next turn, the cops cannot occupy $x_i$ and so cannot surround $R$.

Next suppose that $N[x_i]$ contains a cop for each $i \in \{1, \ldots, \deg(R)\}$, but there is no cop on vertex $R$.  Thus, for each $i \in \{1, \ldots, \deg(R)\}$, $S_i \cup \{x_i\}$ contains a cop.  Note that since there are no vertices in common between $S_i \cup \{x_i\}$ and $S_j \cup \{x_j\}$ for $i \neq j$, each set $S_i \cup \{x_i\}$ contains exactly one cop.  Moreover, each cop is on a vertex of one of the sets $S_i \cup \{x_i\}$, so no vertex in a set $T_j$ contains a cop.
Now, as the robber is not currently surrounded, there must be at least one neighbour of $R$ which is not occupied by a cop; without loss of generality, say $x_1$ has no cop.  This means that there must be a cop on a vertex $s$ of $S_1$.  But since $\delta \geq 3$, $S_1$ has at least one vertex $s' \neq s$, and since there are no cops on $x_1$ or $T_1$, no neighbour of $s'$ contains a cop.  Thus, the robber may move to $x_1$, and will not be surrounded on the next turn as no cop is able to move to $s'$ on the next turn.

Finally, suppose a cop moves to $R$.  Since the robber is not currently surrounded, at least one neighbour of $R$, say $x_1$, has no cop.  The robber will be prevented from moving to $x_1$ only if this causes him to be surrounded immediately or after the cops' next move.  Since $R$ is not adjacent to any element of $S_1$ (else $G$ would contain a 3-cycle), this can only happen if there are $\deg(x_1)-1 \geq \delta-1$ cops in $S_1 \cup T_1$.  But recalling that $S_1 \cup T_1$ is disjoint from $\cup_{i=2}^{\deg(R)} (\{x_i\} \cup S_i \cup T_i)$, it follows that $S_2 \cup T_2$ has no cop.  Thus either $x_1$ or $x_2$ is a vertex to which the robber can move and avoid being surrounded after the cops' next move.
\end{proof}

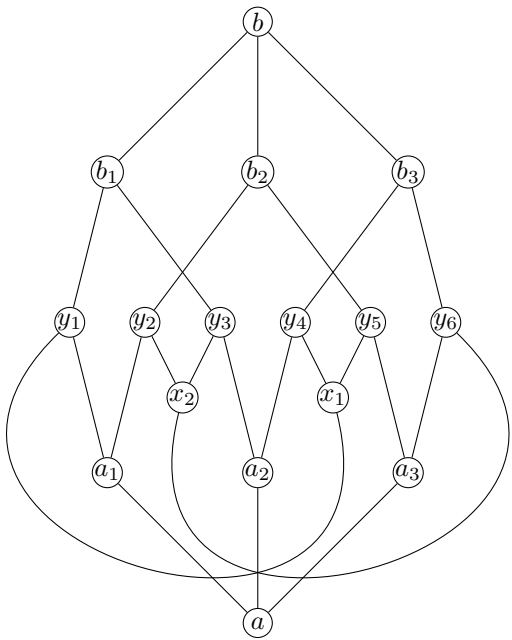
\begin{figure}[t]
\begin{center}
\noindent
\begin{tikzpicture}
		\node [style=vertex] (0) at (0, -3) {{\footnotesize $a$} };
		\node [style=vertex] (1) at (0, 5) {{\footnotesize $b$}};
		\node [style=vertex] (2) at (-2, 3) {{\footnotesize $b_1$ }};
		\node [style=vertex] (3) at (0, 3) {{\footnotesize $b_2$ }};
		\node [style=vertex] (4) at (2, 3) {{\footnotesize $b_3$ }};
		\node [style=vertex] (5) at (-2, -1) {{\footnotesize $a_1$}};
		\node [style=vertex] (6) at (0, -1) {{\footnotesize $a_2$}};
		\node [style=vertex] (7) at (2, -1) {{\footnotesize $a_3$}};
		\node [style=vertex] (12) at (-2.5, 1) {{\footnotesize $y_{1}$}};
		\node [style=vertex] (13) at (-1.5, 1) {{\footnotesize $y_{2}$}};
		\node [style=vertex] (14) at (-0.5, 1) {{\footnotesize $y_{3}$}};
		\node [style=vertex] (15) at (0.5, 1) {{\footnotesize $y_{4}$}};
		\node [style=vertex] (16) at (1.5, 1) {{\footnotesize $y_{5}$}};
		\node [style=vertex] (17) at (2.5, 1) {{\footnotesize $y_{6}$}};
		\node [style=vertex] (18) at (1, 0) {{\footnotesize $x_1$}};
		\node [style=vertex] (19) at (-1, 0) {{\footnotesize $x_2$}};
		\draw (5) to (0);
		\draw (0) to (6);
		\draw (0) to (7);
		\draw (1) to (2);
		\draw (1) to (3);
		\draw (1) to (4);
		\draw (12) to (5);
		\draw (5) to (13);
		\draw (6) to (14);
		\draw (6) to (15);
		\draw (7) to (16);
		\draw (7) to (17);
		\draw (2) to (12);
		\draw (2) to (14);
		\draw (3) to (13);
		\draw (3) to (16);
		\draw (15) to (4);
		\draw (4) to (17);
		\draw (13) to (19);
		\draw (19) to (14);
		\draw [bend right=120, looseness=2.75] (19) to (17);
		\draw (15) to (18);
		\draw (18) to (16);
		\draw [bend right=120, looseness=2.75] (12) to (18);
\end{tikzpicture}
\caption{A cubic graph with girth 6 and surrounding cop number of 3.}\label{fig 1}
\end{center}
\end{figure}

We note that the result of Theorem~\ref{Theorem-Girth} is not true if the hypothesis that $g\geq 7$ is replaced by $g \geq 6$, as we claim the graph $G$ shown in Figure~\ref{fig 1} has girth 6 and $\sigma(P)=\delta(P)=3$.  To see that $G$ has girth 6, note that $G$ is bipartite and clearly has no 4-cycles.  Hence we need to justify the claim that three cops are sufficient to surround the robber on $G$.  One cop (${\cal C}_1$) will play on the vertices of $N[x_1]$ as follows. The cop will initially be placed on $x_1$.  If the robber  ever moves to a vertex in $\{a_1, a_2, a_3, b_1, b_2, b_3\}$, ${\cal C}_1$ will move to the unique neighbour of the robber's position in $N[x_1]$.  For all other robber positions, ${\cal C}_1$ will remain on or move to $x_1$.  Similarly a second cop, say ${\cal C}_2$, will play on the vertices of $N[x_2]$. This cop will initially be placed on $x_2$.  If the robber  ever moves to a vertex in $\{a_1, a_2, a_3, b_1, b_2, b_3\}$, ${\cal C}_2$ will move to the unique neighbour of the robber's position in $N[x_2]$ and, for all other robber positions, ${\cal C}_2$ will remain on or move to $x_2$. The movements of cops ${\cal C}_1$ and ${\cal C}_2$ ensure that once the robber is in $N[a]$ or $N[b]$, he cannot leave. The third cop, ${\cal C}_3$, will move towards the robber's position, to ensure (if necessary) that the robber cannot indefinitely stay on a vertex in $\{y_1, y_2, y_3, y_4, y_5, y_6\}$. Once the robber is in $N[a]$ or $N[b]$,  ${\cal C}_3$ moves to $a$ or $b$ respectively, and the robber is surrounded by the next turn.

With the exception of Lemma~\ref{Lemma-UpperBound},
each of the bounds on $\sigma(G)$ that we have so far seen has been a lower bound.
In Theorem~\ref{Thm-Treewidth}, we will show that
the treewidth of $G$ can be used to provide an upper bound on the surrounding cop number.  
We first review this parameter, previously introduced by Robertson and Seymour~\cite{RS1986}.

Given a graph $G = (V, E)$, a {\em tree decomposition} is a pair $(X, T)$, where $X = \{X_1, \ldots, X_n\}$ is a family of subsets of $V$
called \emph{bags}, and $T$ is a tree whose vertices are the subsets $X_i$, satisfying the following three properties:

\begin{enumerate}
\item $V=\bigcup_{i=1}^nX_i$; that is, each graph vertex is associated with at least one tree vertex;

\item For every edge $\{v, w\}$ in the graph $G$, there is a bag $X_i$ that contains both $v$ and $w$;

\item If $X_i$, $X_j$ and $X_k$ are bags, and $X_k$ is on the path in $T$ from $X_i$ to $X_j,$ then $X_i\cap X_j \subseteq X_k$.
\end{enumerate}

The {\em width} of a tree decomposition is one less than the size of its largest bag $X_i$. The {\em treewidth} of a graph $G,$ written
$\tw(G),$ is the minimum width among all possible tree decompositions of $G$. Note that, by condition (3), for each vertex in $G$, the set of bags
containing the vertex is a subtree of $T$.
It then follows from condition (2) that vertices are adjacent in $G$ only when their corresponding subtrees have a tree vertex in common.

An upper bound involving treewidth is known for the cop number $c(G)$.

\begin{theorem}[\cite{JKT2010}]
For any graph $G$, $c(G) \leq \frac{\tw(G)}{2}+1$.
\end{theorem}

We now establish an upper bound on the surrounding cop number.

\begin{theorem}\label{Thm-Treewidth}
For any graph $G$, $\sigma(G) \leq \tw(G)+1$.
\end{theorem}

\begin{proof}
Let $T$ be a fixed tree decomposition of $G$ such that $T$ has width $\tw(G)$. Place a cop on each vertex of an arbitrary bag $B$ of $T$ and any
remaining cops on arbitrary vertices of $G$.
Hence the robber cannot start on a vertex of $B$. Let $B'$ be the unique bag adjacent to $B$ in $T$ which is on a shortest path
connecting the bag $B$ to the subtree $T'$ of $T - B$ that contains the robber.

The cops now perform a sequence of moves that culminate with a cop being present on each vertex of $B'$.
They do so in such a way that throughout this process a cop is always present on each vertex of $B\cap B'$.
As each edge of $G$ is in a bag, $B\cap B'$ is a cut set of $G$.
At the end of this sequence of moves, when the vertices of $B'$ are occupied by cops,
the cops have therefore restricted the robber's possible positions to the vertices of the set $\bigcup_{X \in T'} X \setminus B'$.

The cops then similarly move to occupy each vertex of $B''$, the unique bag adjacent to $B'$ which is on a shortest path in $T$ to the subtree $T''$ of $T-B'$
that contains the robber, in such a way that throughout this process a cop is always present on each vertex in $B'\cap B''$.
By iterating this process, the cops are able to further restrict the possible positions of the robber.
As $G$ is finite, the robber is eventually surrounded.
\end{proof}

\begin{corollary}\label{Cor-chordal}
If $G$ is a chordal graph then $\sigma(G) \in \{\omega(G)-1,\omega(G)\}$.
\end{corollary}

\begin{proof}
If $G$ is a chordal graph then $\tw(G) = \omega(G)-1$.
It follows from Lemma~\ref{Lemma-OmegaMinusOne} and Theorem~\ref{Thm-Treewidth} that
$\omega(G)-1 \leq \sigma(G) \leq \tw(G) + 1 = \omega(G)$.
\end{proof}

Observe that each of the two possibilities in Corollary~\ref{Cor-chordal} can be realised
when considering trees, which are chordal graphs and have treewidth 1.  If $G$ is a star then $\sigma(G)=1=\omega(G)-1$,
whereas if $G$ is a tree which is not a star then $\sigma(G)=2=\omega(G)$.

Also note that the treewidth of a graph $G$ does not provide a lower bound on $\sigma(G)$ in general.
For instance, observe that $\tw(P_n \Box P_n)=n$, yet $\sigma(P_n \Box P_n) = 3$ when $n \geq 4$
(see Theorem~\ref{Thm-CartesianProductPath}).

We conclude this section with an example of a class of graphs for which the cop number is bounded by a constant,
but the surrounding cop number is on the order of the square root of the number of vertices of the graph.
For a graph $G$, let $L(G)$ denote the {\em line graph} of $G$, namely the graph having vertex set $V(L(G)) = E(G)$
and such that two vertices of $L(G)$ are adjacent in $L(G)$ if (as edges) they share a vertex in $G$.

\begin{theorem}\label{Thm-LineGraphKn}
For any integer $n\geq 3$, $\sigma(L(K_n))=2(n-2)$.
\end{theorem}

\begin{proof}
It is clear from Lemma~\ref{Lemma-deltabound} that $\sigma(L(K_n))\geq 2(n-2)$ since $L(K_n)$ is a $2(n-2)$-regular graph.
It is easily verified that $\sigma(L(K_3)) = 2$.
Hence it remains to show that $2(n-2)$ cops suffice when $n\geq 4$.
Let $V(K_n)=\{1,\ldots,n\}$, so  that $V(L(K_n))= \big\{ \{x,y\} : x,y\in V(K_n), x\neq y \big\}$.

Initially place one cop on each of the $2(n-2)$ vertices of
$\big\{ \{1,x\} : 2 \leq x \leq n \big\} \cup \big\{ \{2,x\} : 3 \leq x \leq n-1\} \big\}$.
Then without loss of generality we can assume the robber starts on vertex $\{2,n\}$, $\{3,n\}$ or $\{3,4\}$.
We will show that in each case the cops can immediately move to surround the robber.

If the robber starts on vertex $\{2,n\}$ then for each $x \in \{3,4,\ldots,n-1\}$ the cop on vertex $\{1,x\}$ moves to $\{x,n\}$,
causing the robber to become surrounded.
Now suppose the robber starts on vertex $\{3,n\}$.
For each $x \in \{4,5,\ldots,$ $n-1\}$ the cop on vertex $\{2,x\}$ moves to $\{3,x\}$.
Simultaneously, for each $y \in \{2,3,\ldots,n-1\} \setminus \{3\}$ the cop on vertex $\{1,y\}$ moves to $\{y,n\}$.
The robber is then surrounded.
Finally, suppose the robber starts on vertex $\{3,4\}$.
For each $x \in \{5, 6, \ldots,$ $n-1\}$ the cop on vertex $\{1,x\}$ moves to $\{3,x\}$ and the cop on $\{2,x\}$ moves to $\{4,x\}$.
Simultaneously the cop on $\{1,n\}$ moves to $\{3,n\}$,
the cop on $\{1,2\}$ moves to $\{1,4\}$ and the cop that was initially on $\{1,4\}$ moves to $\{4,n\}$.
The robber is then surrounded.
\end{proof}

Note that the line graph of $K_n$ has order $\frac{n(n-1)}{2}$ and so the surrounding cop number of $L(K_n)$
is $\Theta\big(\sqrt{|V(L(K_n))|}\,\big)$,
whereas the cop number of $L(K_n)$ is bounded by a constant.
Specifically,  $c(L(K_n))\leq 2$ because $c(L(G))\leq c(G)+1$ for any graph $G$~\cite{Dudek2014}.

\section{Algorithmic considerations} \label{Section-Algorithmic}

In~\cite{BC2009}, Bonato and Chiniforooshan showed that a robber can win the original Cops and Robbers game when $k$ cops are deployed
if and only if there exists a function $\psi$ from the set of all possible placements of cops to
the power set $\mathcal{P}(V(G))$ such that $\psi$ has certain properties, one of which is that $\psi(T) \neq \emptyset$ for each configuration $T$.
More technically, the set of all possible placements of $k$ cops within $G$ can be represented as the vertex set of the graph
$\boxtimes^k G = G \boxtimes G \boxtimes \cdots \boxtimes G$ (with $k-1$ instances of the strong product operation $\boxtimes$
being applied), where two vertices are adjacent if the cops can move between the two configurations in a single move.
For each $T$, the set $\psi(T)$ ultimately represents those vertices of $G$ that the robber ought to move onto when the cops are at $T$.

In this section we adapt results from~\cite{BC2009} to the new context of the surrounding cop number.
We begin with a characterisation of when the robber can win when in the presence of $k$ cops.

\begin{theorem}
\label{Theorem-BC2009adapted}
Let $k \in \mathbb{N}$. Then $\sigma(G) > k$ if and only if there exists a function $\psi : V(\boxtimes ^k G) \rightarrow \mathcal{P}(V(G))$ with the following three properties:

\begin{enumerate}

\item[(1)] $\forall T \in V(\boxtimes ^k G)$, $\psi (T) \neq \emptyset$;

\item[(2)] $\forall T \in V(\boxtimes ^k G)$, $\psi (T) \subseteq V(G)\setminus (A_T \cup B_T \cup C_T)$ where $A_T$ is the set of all vertices of $G$ that are in the $k$-tuple $T$, $B_T$ is the set of vertices in $G$ that are currently surrounded by $T$, and $C_T$ is the set of vertices in $G$ that can be surrounded by $T$ in one move;

\item[(3)] $\forall T_1 T_2 \in E(\boxtimes ^k G)$, $\psi (T_1) \subseteq N_G[\psi (T_2)]$.

\end{enumerate}
\end{theorem}

In the interest of brevity we omit a proof, which can be obtained by suitably adapting the proof found in~\cite{BC2009}
concerning the cop number.

\begin{corollary}
\label{Corollary-psiempty}  Let $G$ be a connected graph, and let $\psi: V(\boxtimes ^k G) \rightarrow \mathcal{P}(V(G))$ be a function satisfying properties (2) and (3) of Theorem~\ref{Theorem-BC2009adapted}.
If there exists $T_1 \in V(\boxtimes ^k G)$ with $\psi (T_1)=\emptyset$,
then $\psi (T) =\emptyset$ for each $T \in V(\boxtimes ^k G)$.
\end{corollary}

\begin{proof}
Since $G$ is connected, so is $\boxtimes ^k G$. Thus there is a finite walk $W={T_1, T_2,...,T_m}$ in $\boxtimes ^k G$ which starts at $T_1$ and visits every vertex of $\boxtimes ^k G$ at least once.  Suppose that for some $i \in \{1, \ldots, m\}$, $\psi(T_i)=\emptyset$.  Noting that $T_i T_{i+1} \in E(\boxtimes ^k G)$, property (3) implies that $\psi(T_{i+1}) \subset N_G[\psi(T_i)] = \emptyset$, and thus $\psi(T_{i+1}) = \emptyset$.  But recalling that $\psi(T_1) = \emptyset$, it follows by induction that $\psi(T_i) = \emptyset$ for all $i \in \{1, \ldots, m\}$.  Since $W$ contains each vertex of $\boxtimes ^k G$, the result follows.
\end{proof}

Bonato and Chiniforooshan~\cite{BC2009} also presented an algorithm to decide whether $k$ cops suffice to capture a robber in a given graph.
With only minor alterations to the way in which the algorithm is initialised,
it can also be used to determine the surrounding cop number $\sigma(G)$ of a graph $G$
(see Algorithm~\ref{alg1}).

\vspace{0.5\baselineskip}
\begin{algorithm}[H]
\label{alg1}

\SetKwInOut{Input}{input}\SetKwInOut{Output}{output}

\Input{connected graph $G=(V,E)$, number of cops $k \in \mathbb{N}$}

\Output {either $\sigma(G) \leq k$ or $\sigma(G)>k$}

\BlankLine

initialise $\psi (T)$ to $V(G)\setminus (A_T \cup B_T \cup C_T)$ for all $T \in V(\boxtimes ^k G)$

\Repeat{the value of $\psi$ is unchanged}{

\ForAll{$T T^\prime \in E(\boxtimes ^k G)$}{

$\psi (T) \leftarrow \psi (T) \cap N_G [\psi (T^\prime)]$

$\psi (T^\prime) \leftarrow \psi (T^\prime) \cap N_G [\psi (T)]$

}

}

\uIf{$\exists T \in V(\boxtimes ^k G)$ such that $\psi (T) = \emptyset$}{

\Return{$\sigma(G) \leq k$} \
\text{(cops win)}

}

\Else{

\Return{$\sigma(G) > k$} \
\text{(robber wins)}

}

\caption{check potential $\sigma(G)$ value}
\end{algorithm}

\begin{theorem}
\label{Theorem-algruntime}
Algorithm 1 runs in time $\bigO(kn^{k+3}+n^{3k+3})$ where $n=|V(G)|$.
\end{theorem}

We omit the proof of Theorem~\ref{Theorem-algruntime} as it follows in a similar fashion to the proof of Theorem 2.2 in~\cite{BC2009}.

We implemented Algorithm~\ref{alg1}, as well as the corresponding algorithm from~\cite{BC2009} for the cop number.
When testing to see whether $k$ cops are able to win, the algorithm performs best when $k$ is small.
For other values of $k$, the order of the graph $\boxtimes^k G$ makes it impractical to execute the algorithm.
Nevertheless, we found some interesting empirical results, particularly for generalised Petersen graphs,
which we will discuss in more detail in Section~\ref{Section-GeneralisedPetersen}.

We conclude this section with some comments about the computational complexity
of determining the surrounding cop number, expressed as a decision problem as follows.

\begin{tabbing}
\hspace*{9mm}\=
{\textsc{Surrounding Cop Number}}\\
\>\hspace*{5mm}\=
{Instance}: A graph $G$ and an integer $k$.\\
\>\>
{Question}: \= Is $\sigma(G) \leq k$?
\end{tabbing}

It was shown in 2010 that
the analogous
\textsc{Cop Number} decision problem for
determining the cop number $c(G)$ of an arbitrary graph $G$ is {\NP}-hard~\cite{FGKNS2010},
although whether the problem is in {\NP} remains an open question.
It was conjectured in~\cite{GR1995} that the problem of determining $c(G)$ is {\EXPTIME}-complete; the truth of this conjecture was established by Kinnersley~\cite{Kinnersley2015}.

Theorem~\ref{Theorem-algruntime} shows that for fixed values of $k$
the \textsc{Surrounding Cop Number} problem can solved in polynomial time.
However, for arbitrary $k$ the computational complexity of the
problem is not currently known.
We leave as open problems the questions of determining whether
\textsc{Surrounding Cop Number}
is in {\NP}, is {\NP}-hard, or is {\EXPTIME}-complete.

\section{Classes of Graphs}
\label{Section-GraphClasses}

It is natural to ask how the surrounding cop number behaves for various classes of graphs, as well as in relation to some common graph operations.
For instance, given two graphs $G$ and $H$, we consider the surrounding cop numbers of product graphs built from $G$ and $H$.
We also consider graphs arising from designs, as well as generalised Petersen graphs.

\subsection{Graph Products} \label{Section-Product}
We recommend the monograph by Imrich and Klav\v{z}ar~\cite{ImrichKlavzar} for an overview of graph products and for the notation that we follow.
Here we consider
the Cartesian product $G \Box H$,
the strong product $G \boxtimes H$, and the lexicographic product $G \circ H$
for various graphs $G$ and $H$.
We begin with the Cartesian product, for which the following result is known for the cop number.

\begin{theorem}[\cite{Tosic1988}]
If $G$ and $H$ are each connected, then $c(G \Box H) \leq c(G) + c(H)$.
\end{theorem}

We establish a similar result for the surrounding cop number.

\begin{theorem}
\label{Thm-CartesianProduct}
If $G$ and $H$ are each connected, then $\sigma(G \Box H) \leq \sigma(G) + \sigma(H)$.
\end{theorem}

\begin{proof}
The cops adopt a strategy whereby $\sigma(G)$ cops first seek to be within the robber's copy of $G$,
while $\sigma(H)$ other cops seek to be within the robber's copy of $H$.
Note that
any time that the robber stays within his present copy of $G$,
the $\sigma(G)$ $G$-seeking cops can each move along an edge within their copy of $H$ and thereby get closer to the robber's copy of $G$,
and similarly the $\sigma(H)$ $H$-seeking cops will get closer to the robber's copy of $H$ whenever the robber stays within his present copy of $H$.
Clearly at least one of the two teams of cops will achieve their initial goal after a finite number of rounds of the game.

Without loss of generality (due to the symmetry between $G$ and $H$),
we may now suppose that the $\sigma(G)$ cops that sought to occupy the same copy of $G$ as the robber have achieved their initial goal.
Henceforth these $\sigma(G)$ cops mirror any move that the robber makes when he moves to a different copy of $G$, or else they follow a strategy within
their copy of $G$ that leads to them occupying each of the robber's neighbouring vertices within that copy of $G$
(which, by definition, can be achieved by $\sigma(G)$ cops).
Recall that any time that the robber either stays still or moves within his copy of $H$, the $H$-seeking cops are able to get closer to the robber's copy of $H$.
Observe that the strategy now employed by the $G$-seeking cops will therefore enable the $H$-seeking cops to successfully reach the robber's copy of $H$ in a finite number of moves,
at which time both teams of cops will have achieved their initial goal.

Both teams now follow similar strategies; that is, if the robber moves to a new copy of $G$ (resp.~$H$) then the $G$-seeking (resp.~$H$-seeking) cops
make a parallel move, while the $H$-seeking (resp.~$G$-seeking) cops work towards surrounding the robber within his copy of $H$ (resp.~$G$).
Clearly one of these teams, say the $G$-seeking team, will achieve the goal of surrounding the cop within a copy of $G$
after a finite number of moves.  Once that has happened, the robber's movement is restricted to a single copy of $H$,
and the team of $\sigma(H)$ cops can then follow a strategy to ultimately surround the robber within that copy of $H$.
\end{proof}

We now show that the bound of Theorem~\ref{Thm-CartesianProduct} is sometimes, but not always, an equality.
For instance, equality is obtained when $G=H=P_2$, but not when $G$ and $H$ are paths on four or more vertices.

\begin{theorem} \label{Thm-CartesianProductPath}
Let $2 \leq m \leq n$ be integers.
If $m,n \leq 3$ then $\sigma(P_m \Box P_n) = 2$;
otherwise $\sigma(P_m \Box P_n) = 3$.
\end{theorem}

\begin{proof}
It is easy to verify that $\sigma(P_m \Box P_n) = 2$ when $m,n \leq 3$.
So we henceforth assume $n \geq 4$.
We commence by showing that two cops are unable to surround the robber.
As a winning strategy for the robber, he initially chooses a position that is not one of the four vertices of degree 2.
The cops cannot surround the robber on such a vertex, although they can force the robber to move to a new vertex.
When the robber is forced to move, one cop shares the robber's vertex, which has at least two unoccupied neighbouring vertices, at most one of which has degree 2.
Thus the robber is able to move to a new vertex that is not a vertex of degree 2.

We now present a strategy that enables three cops to surround the robber.
In reference to the $m \times n$ rectangular grid of vertices formed by the graph $P_m \Box P_n$,
for their initial positions,
cop ${\cal C}_1$ positions herself at the rightmost vertex of the top row,
cop ${\cal C}_2$ is immediately below ${\cal C}_1$,
and cop ${\cal C}_3$ is at the leftmost vertex of the top row.
Initially the cops ${\cal C}_1$ and ${\cal C}_2$ sweep leftwards until they are in the same column of vertices as the robber, while ${\cal C}_3$ remains stationary.

Once ${\cal C}_1$ and ${\cal C}_2$ are in the same column as the robber, the cops adjust their strategy.
First, suppose that the robber is on the topmost row when ${\cal C}_1$ and ${\cal C}_2$ move to occupy his column; thus, ${\cal C}_1$ will occupy the robber's vertex and ${\cal C}_2$ will be directly below him, forcing the robber to move horizontally.  In subsequent rounds, ${\cal C}_1$ and ${\cal C}_2$ move horizontally to follow the robber to his new column; for the rest of the game, the movements of ${\cal C}_1$ and ${\cal C}_2$ contain the robber to the topmost row and prevent him from passing on any turn.  In each subsequent turn, ${\cal C}_3$ moves rightwards on the top row, unless such a move would cause her to occupy the same vertex as the robber, in which case she remains stationary.  Eventually, ${\cal C}_3$ will occupy the robber's left neighbour, and the robber will thereafter be forced to move rightwards.  When he reaches the upper right corner, he will be surrounded on the cops' next turn.

We now assume that the robber is not on the topmost row when ${\cal C}_1$ and ${\cal C}_2$ move to occupy his column.  In the next phase of the game, cops ${\cal C}_1$ and ${\cal C}_2$ move horizontally to follow the robber to a new column if he chooses to move to a new column; otherwise they do not move.
Meanwhile, ${\cal C}_3$ moves down from the topmost row until she finds herself in the same row as the robber.

Before the robber makes his next move, he finds himself in the same row as ${\cal C}_3$,
and in the same column as ${\cal C}_1$ and ${\cal C}_2$. 
Note that ${\cal C}_3$ either occupies the same vertex as the robber in the leftmost column or is to the left of the robber, while ${\cal C}_1$ and ${\cal C}_2$ are either both above the robber, or else ${\cal C}_1$ is his above neighbour and ${\cal C}_2$ occupies the same vertex as the robber.
If the robber moves vertically, ${\cal C}_3$ then makes a parallel move to the robber's new row, while
${\cal C}_1$ and ${\cal C}_2$ both move down.
If the robber moves horizontally, then on their next move ${\cal C}_1$ and ${\cal C}_2$ move horizontally to the robber's new column,
while ${\cal C}_3$ either moves to the right if the vertex to her right is not occupied by the robber,
or else ${\cal C}_3$ remains stationary as the robber's left neighbour.
If the robber should happen to remain stationary, then ${\cal C}_1$ and ${\cal C}_2$ both move down,
and ${\cal C}_3$ moves to the right if the robber is not on the vertex that is her right neighbour.
By moving in this fashion, the cops will guide the robber to move downward and to the right,
ultimately surrounding him in the rightmost vertex of the bottom row.
\end{proof}

On the matter of the strong product, we first show that exact values for the surrounding cop number are able to be determined
when taking the product of two paths.

\begin{theorem}\label{Thm-StrongProduct}
Let $2 \leq m \leq n$ be integers.  Then
\[
\sigma(P_m \boxtimes P_n) = \left \{ \begin{array}{ll}
5, & \mbox{if } m \geq 4 \\
4, & \mbox{if } m=3, \mbox{ or } m=2 \mbox{ and } n \geq 4 \\
3, & \mbox{if } m=2 \mbox{ and } n \leq 3.
\end{array}
\right.
\]
\end{theorem}

\begin{proof}
For convenience, label the vertices of $P_m$ (resp.\ $P_n$) with integers from 1 to $m$ (resp.\ 1 to $n$) so that the vertices of $P_m \boxtimes P_n$ are labelled with 2-dimensional coordinates from the set $\{1,2,\ldots,m\} \times \{1,2,\ldots,n\}$ with $(1,1)$ at the bottom left.  It is easy to verify that $\sigma(P_m \boxtimes P_n) = 3$ when $m=2$ and $n \leq 3$, so we henceforth assume that if $m=2$, then $n \geq 4$.

Suppose now that $m=2$ and $n \geq 4$, or that $m=3$.
If only three cops are available, then the robber can initially position himself on a vertex that is not one of the four vertices of degree 3.
The robber cannot be surrounded by three cops when he is on a vertex of degree exceeding 3
and, when forced to move from his current position, at least one neighbouring vertex of degree greater than 3 is not occupied by a cop and so the robber
moves to such a vertex.  The robber can iterate this manoeuvre indefinitely to win the game.

If $m=2$ and $n \geq 4$ then a winning strategy for four cops is to initially place one cop on each vertex of the two bottommost columns (that is, on vertices $(1,1)$, $(2,1)$, $(1,2)$ and $(2,2)$), 
and at each turn they move upwards, ultimately corralling the robber into the topmost copy of $P_2$ where he will be easily surrounded.

If $m=3$ then a winning strategy for four cops is as follows.  Three of the cops, ${\cal C}_1, {\cal C}_2$ and ${\cal C}_3$, initially place themselves on vertices $(1,1)$, $(2,1)$ and $(3,1)$, respectively, and the fourth cop, ${\cal C}_4$, begins on the vertex $(2,2)$.  The cops' first aim is to force the robber onto the topmost row.  Initially, on each turn each cop moves one vertex upward, until ${\cal C}_4$ occupies the same row as the robber.  Say this occurs when the robber is on vertex $(i,j+1)$, so that ${\cal C}_1, {\cal C}_2, {\cal C}_3$ are on vertices $(1,j), (2,j)$ and $(3,j)$, respectively, and ${\cal C}_4$ is on vertex $(2,j+1)$.  If the robber moves to a vertex of row $j+2$ on his next move, the cops continue their upward trajectory.  Otherwise, after his next move the robber is on either $(1,j+1)$ or $(3,j+1)$.  By symmetry we may assume the robber is on $(1,j+1)$.  Now, ${\cal C}_3$ moves to $(2,j+1)$ and ${\cal C}_4$ moves to $(1,j+1)$.  The robber is now forced to move into column $j+2$.  The cops ${\cal C}_1, {\cal C}_2, {\cal C}_3$ and ${\cal C}_4$ move to $(1,j+1), (2,j+1), (3,j+1)$ and $(2,j+2)$, respectively, thus resetting their formation.  In either case, they have forced the robber to move further upwards.

Eventually, after a cops' turn, the robber will occupy a vertex of the form $(i,n)$, and the cops will occupy vertices $(1,n-1), (2,n-1), (3,n-1)$ and $(2,n)$.  The robber is either already surrounded, or else he is forced to move to a corner vertex where he will be surrounded.

Finally assume that $4 \leq m \leq n$.
To show that five cops suffice, we present a strategy that shows how they can always surround a robber.
First position the five cops  $\cC_1,\ldots,\cC_5$ on vertices $(m,n),(m,n-1),(m,n-2),(m-1,n)$ and $(m-1,n-1)$ respectively. For most of their strategy, the cops will maintain this configuration relative to each other.  The cops' initial formation is illustrated in Figure~\ref{Figure:StrongGrid}.
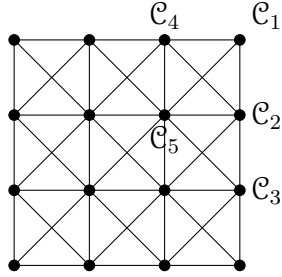
\begin{figure}[ht]
\begin{center}
\begin{tikzpicture}[x=1cm,y=1cm,scale=1]
\foreach \x in {0,...,3} \foreach \y in {0,...,3}{
	\draw[fill=black](\x,\y) circle (2pt);
}
\foreach \x in {0,...,2} \foreach \y in {0,...,3}{
	\draw (\x,\y) -- (\x+1,\y);
}
\foreach \x in {0,...,3} \foreach \y in {0,...,2}{
	\draw (\x,\y) -- (\x,\y+1);
}
\foreach \x in {0,...,2} \foreach \y in {0,...,2}{
	\draw (\x,\y) -- (\x+1,\y+1);
	\draw (\x,\y+1) -- (\x+1,\y);
}
\draw (2,3) node[above]{${\cal C}_4$};
\draw (3,3) node[above right]{${\cal C}_1$};
\draw (3,2) node[right]{${\cal C}_2$};
\draw (3,1) node[right]{${\cal C}_3$};
\draw (2,2) node[below]{${\cal C}_5$};
\end{tikzpicture}
\caption{Initial positions of cops in $P_m \boxtimes P_n$, $4 \leq m \leq n$.}
\label{Figure:StrongGrid}
\end{center}
\end{figure}
We claim that while maintaining this formation as they move,
the cops are able to achieve one of the following intermediate goals:
\pagebreak
\begin{enumerate}
\item[(i)]  the robber is on vertex $(x,1)$ for some $x \in \{1,2,\ldots,m-1\}$
and the cops are in formation with $\cC_3$ at $(x+1,1)$;
or
\item[(ii)] the robber is on vertex $(1,n)$ and has been surrounded by $\cC_1$, $\cC_2$ and $\cC_5$.
\end{enumerate}

One facet of the cops' strategy while they move in formation
is to keep the robber confined to the rectangular region bounded by $(1,1)$ on the bottom left and the position of $\cC_1$ on the top right
(with the perimeter being included as part of the region described); note that the boundaries of this region will change as the cops move.
Maintaining this restriction is easily accomplished:
if the robber is in the same row as $\cC_1$ then the cops move leftward and upward (unless they are unable to move in one of those directions),
and
if the robber is in the same column as $\cC_1$ then the cops move downward and to the right (unless they are unable to move in one of those directions).
Otherwise, if condition (i) is not met and the robber is in neither the same row nor the same column as $\cC_1$
then the cops move down and/or left toward the robber, reducing both their vertical and horizontal distance between the robber and $\cC_5$ when possible.
Should $\cC_5$ land on the robber, then he will be forced to move down and/or left.

Whenever condition (i) is achieved the cops change strategy, acting to corral the robber to move left towards being ultimately surrounded at $(1,1)$.

All that now remains is to show that four cops are insufficient to surround the robber.
If the robber adopts a strategy whereby he initially positions himself on a vertex that is not one of the four vertices of degree 3,
then he will be on a vertex of degree at least 5.  With four cops, clearly he cannot be surrounded while he is on such a vertex.
However, the cops can force him to move by having a cop occupy his vertex, at which time at least two of his neighbouring vertices
are unoccupied, at least one of which is a vertex of degree 5 or more.
Hence the robber can perpetually move among vertices of degree at least 5 and thereby avoid ever being surrounded.
\end{proof}

Note that
Neufeld and Nowakowski proved the following upper bound for the cop number of the strong product of graphs.

\begin{theorem}[\cite{NeufeldNowakowski1998}]
If $G$ and $H$ are each connected, then $c(G \boxtimes H) \leq c(G) + c(H) - 1$.
\end{theorem}

The following result, however, establishes that a similar upper bound on $\sigma(G \boxtimes H)$ in terms of $\sigma(G)$ and $\sigma(H)$ does not exist.  In particular, although the surrounding cop number of a star is 1, the surrounding cop number of the strong product of a star with itself is bounded below by $n+1$.
Note, moreover, that Theorem~\ref{Thm-StrongProductStar} gives a family of graphs showing that the surrounding cop number of $G$ can be arbitrarily larger than $\omega(G)-1$, as the clique number of $K_{1,n} \boxtimes K_{1,n}$ can easily be shown to be 4.

\begin{theorem} \label{Thm-StrongProductStar}
For $n \geq 1$, $\sigma(K_{1,n} \boxtimes K_{1,n}) \geq n+1$.
\end{theorem}

\begin{proof}
When $n \leq 2$ the result follows from Theorem~\ref{Thm-StrongProduct}, and so we henceforth assume that $n \geq 3$.
Let $v_0$ be the vertex of degree $n$ in $K_{1,n}$, and $L$ be the set of $n$ leaves of $K_{1,n}$. In $K_{1,n} \boxtimes K_{1,n}$, let $A$ be the set of vertices of the form $(v_0,\ell)$ with $\ell\in L$ and $B$ be the set of vertices of the form $(\ell,v_0)$ with $\ell\in L$. Observe for every $\ell,\ell'\in L$, $(v_0,\ell)$ is adjacent to $(\ell', v_0)$.  In other words, every vertex in $A$ is adjacent to every vertex in $B$.  As $|A|=|B|=n$, thus $K_{n,n}$ is a subgraph of $K_{1,n}\boxtimes K_{1,n}$ not containing $(v_0,v_0)$. As the robber can just play on this subgraph, at least $n$ cops are needed to surround the robber.  Since $(v_0,v_0)$ is a universal vertex, an additional cop must always stay on $(v_0,v_0)$ (or the robber has an escape).
\end{proof}

This result can be generalised further to the strong product of a complete bipartite graph with itself.
\begin{theorem}
If $a,b \geq 1$ then $\sigma(K_{a,b} \boxtimes K_{a,b}) \geq ab$.
\end{theorem}

\begin{proof}
Let the bipartition of $K_{a,b}$ be $(A,B)$, where $|A|=a$ and $|B|=b$.  Let $X$ be the set of vertices of $K_{a,b} \boxtimes K_{a,b}$ of the form $(\alpha,\beta)$ and $Y$ be the set of vertices of the form $(\beta,\alpha)$, where $\alpha \in A$ and $\beta \in B$.  Note that $|X|=|Y|=ab$.  In $K_{a,b} \boxtimes K_{a,b}$, every vertex in $X$ is adjacent to every vertex of $Y$, giving a subgraph isomorphic to $K_{ab,ab}$.  Considering that the robber can always play on this subgraph, the result follows.
\end{proof}

Finally we consider the lexicographic product of two graphs.

\begin{theorem}\label{Thm-LexicographicProductBound}
If $G$ and $H$ are connected graphs then $\sigma(G \circ H)\leq \sigma(G)|V(H)|+\sigma(H)$.
\end{theorem}

\begin{proof}
Begin by selecting $\sigma(G)$ copies of $H$, say $H_1, H_2,\ldots, H_{\sigma(G)}$, and then place a cop on each vertex of $\bigcup_{i=1}^{\sigma(G)} V(H_i)$.
By imitating a strategy for surrounding a robber in $G$,
trap the robber on one copy of $H$; that is if a cop moves from $x$ to $y$ in $G$, then move all of the cops from $H_x$ to $H_y$ in $G \circ H$.
Subsequently use the remaining $\sigma(H)$ cops to play on the copy of $H$ containing the robber.
\end{proof}

We note that $\delta(G \circ H)=\delta(G) |V(H)|+\delta(H)$, and so the bound given by Theorem~\ref{Thm-LexicographicProductBound} is tight
whenever $\sigma(G)=\delta(G)$ and $\sigma(H)=\delta(H)$.


\subsection{Graphs arising from designs}

In~\cite{BonatoBurgess}, Bonato and Burgess investigated the cop number for several classes of graphs that are based on combinatorial designs.
Accordingly we now proceed to consider the surrounding cop number for some of these classes of graphs.
Before discussing them, though, we introduce several definitions.

For integers $v > k \geq 2$ and $\lambda \geq 1$ we define a {\em balanced incomplete block design} with parameters $(v,k,\lambda)$,
denoted as a BIBD$(v,k,\lambda)$, to be a pair $(X,\mathcal{B})$ where $X$ is a $v$-set and $\mathcal{B}$ is a collection of $k$-subsets of $X$ such that
each 2-subset of $X$ is contained within exactly $\lambda$ of the elements of $\mathcal{B}$.
The elements of $\mathcal{B}$ are called the {\em blocks} of the design.
In any BIBD$(v,k,\lambda)$, each element or {\em point} of $X$ must occur in a constant number of blocks, namely
$\frac{\lambda (v-1)}{k-1}$; this value is known as the design's {\em replication number}, often represented as $r$.
A BIBD$(n^2 + n + 1, n+1, 1)$ is known as a {\em projective plane of order $n$}.
If the blocks of a BIBD$(v,k,\lambda)$ can be partitioned into sets of blocks, each of which contains each point of $X$ exactly once amongst its blocks,
then such a partition is called a {\em resolution} of the design into {\em parallel classes}.
Any BIBD that admits a resolution is a {\em resolvable} BIBD, often denoted as RBIBD.
Given a BIBD$(v,k,\lambda)$, say $(X,\mathcal{B})$, its {\em block-intersection graph} is the graph having $\mathcal{B}$ as its vertex set,
and two vertices are adjacent if (as blocks) they have non-empty intersection.
Another graph that we can construct from a BIBD$(v,k,\lambda)$ $(X,\mathcal{B})$ is its {\em incidence graph}, for which $X \cup \mathcal{B}$ is the vertex set
and a vertex $x \in X$ is adjacent to a vertex $B \in \mathcal{B}$ if and only if $x \in B$.  For further background on design theory, see~\cite{Handbook}.

The cop number of the incidence graph of a projective plane of order $k-1$ (and hence having block size $k$) is known to be $k$;
for a proof see Baird~\cite{Baird2011} or Pra\l{}at~\cite{Pralat2010}.
For the game of Surrounding Cops and Robbers, one additional cop is needed.

\begin{theorem}\label{Thm-ProjectivePlane}
If $G$ is the incidence graph of a projective plane with block size $k$, then $\sigma(G)= k+1$.
\end{theorem}

\begin{proof}
Let $(X,\mathcal{B})$ be a projective plane with block size $k$, and incidence graph $G$.
We first show that $k+1$ cops suffice.
Let $x\in X$ and, noting that there are exactly $k$ blocks containing $x$, put a cop on $x$ and on each $B\in\mathcal{B}$ such that $x\in B$.
If the robber starts on a block $B_0$, then the cops on blocks move to immediately surround the robber because $B_0$ intersects with each of the original blocks at distinct points.

So suppose the robber starts on a point $y\neq x$. As long as the cops occupy all blocks $B$ such that $x\in B$, then the robber stays on $y$ until forced to move (otherwise he will be immediately surrounded as in the previous case).
The cop on  point $x$ moves to $B'$ such that $\{x,y\}\subseteq B'$ and then moves to $y$;  meanwhile, the other cops remain stationary.  The robber is then forced to move to a block containing $y$, and is then surrounded as above.

We now assume towards a contradiction that $\sigma(G)\leq k$. If $k$ cops suffice, then without loss of generality, in the final position the robber is on a point $x$ and the $k$ cops are on each of the blocks containing $x$.  (If the robber is surrounded on a block, we consider the game on the dual.)
If the game ends with the robber's turn then, in the previous move, the robber was on a block $B$ and the cops move to all blocks containing $x$, including $B$.  Note that $x$ is the only point of $B$ which is surrounded, so the robber would have chosen another vertex $y \neq x$ in $B$.

Now suppose the game ends when the robber is surrounded after the cops' turn.
Consider the $k$ blocks $B_1, \ldots, B_k$  which contain $x$.  At the end of the final turn, the $k$ cops must each occupy these blocks.  In the move prior to surrounding $x$, we can assume that $\ell$ cops are on points and $k-\ell$ cops are on blocks, where $1 \leq \ell \leq k$.  Say cops ${\cal C}_1, \ldots, {\cal C}_{\ell}$ occupy points $a_1, a_2, \ldots, a_{\ell}$, respectively, such that $x \neq a_{i} \in B_{i}$, and $k-\ell$ cops occupy the blocks $B_{\ell+1}, \ldots, B_k$.
Note that for each $i \in \{1, \ldots, \ell\}$, no cop other than ${\cal C}_i$ can occupy a point $y \in B_i$, as $B_i$ is the only block containing both $y$ and $x$; moreover, no cop other than ${\cal C}_i$ can be on a block $B'$ which intersects $B_i \setminus\{x\}$.  Thus, if the robber had passed on his previous turn, he could instead have moved to one of the blocks $B_1, \ldots, B_{\ell}$ and would have been safe for another round.  Similarly, had he occupied one of the blocks $B_1, \ldots, B_{\ell}$, he could have passed rather than move to $x$.

Now suppose that the robber had occupied a block $\{B_{\ell+1}, \ldots, B_k\}$, and was forced to move from that block by a cop. Without loss of generality, suppose the robber moved to $x$ from $B_k$.  If $\ell< k-1$, then choose any point $z \in B_k \setminus \{x\}$.  Note that since $B_{\ell+1}$ and $B_k$ share only vertex $x$, $B_{\ell+1}$ and $z$ are not adjacent in the incidence graph.  Thus, the robber could have moved to $z$ rather than $x$; the cop on $B_{\ell+1}$ could not occupy his neighbour on the next round, so he would have avoided capture.

Finally, suppose $\ell=k-1$ and the robber was forced to move from $B_k$ to $x$ by a cop occupying $B_k$.  Thus, after the robber's move to $x$, the other $k-1$ cops occupy vertices $a_1, \ldots, a_{k-1}$ of blocks $B_1, \ldots, B_{k-1}$, respectively.  Consider the block $\hat{B}$ containing points $a_1$ and $a_2$.  Since $(X,\mathcal{B})$ is a projective plane, $\hat{B}$ and $B_k$ must intersect in exactly one point, say $z$.  Noting that $\hat{B} \notin \{B_1, \ldots, B_k\}$ (as $a_1$ and $a_2$ are in distinct blocks containing $x$), we have that $z \neq x$.  But the cops occupying $a_1$ and $a_2$ cannot move to distinct neighbours of $z$; thus, the robber could have moved to $z$ rather than $x$ and remained safe for another round.
\end{proof}

\begin{theorem}
Let $v > k \geq 2$ be integers such that $v\equiv 0$ (mod $k$).
If $G$ is the incidence graph of a resolvable BIBD($v,k,1$),
then $\sigma(G)= \tfrac{v}{k}+1$.
\end{theorem}

\begin{proof}
It can be seen that $\tfrac{v}{k}+ 1$ cops suffice
by placing one cop on each of the $\tfrac{v}{k}$ blocks of a parallel class $\Pi$ and one cop on a point $x$,
and proceeding as in Theorem~\ref{Thm-ProjectivePlane}.

We now show that if there are at most $\frac{v}{k}$ cops, then the robber is always able to return to a point-vertex and hence avoid being surrounded.
Note that point-vertices have degree $\frac{v-1}{k-1}$, and that $\frac{v-1}{k-1}>\frac{v}{k}$ since $v>k$.
The robber adopts a strategy whereby he always stays on (or moves to) a point-vertex, unless forced to move to a block-vertex.
Suppose that at some point the robber is forced to move to a block-vertex because a cop has moved to his current position at point-vertex $x$.
Let $t$ denote the number of neighbouring block-vertices that are occupied by cops, so that the robber has $\frac{v-1}{k-1} - t$ valid moves.
The corresponding $\frac{v-1}{k-1} - t$ blocks collectively contain $(\frac{v-1}{k-1} - t)(k-1)$ points other than point $x$, each of which
has a point-vertex that must be occupiable by a cop if the cops are to be able to surround the robber once he has moved away from $x$.
Moreover, these $(\frac{v-1}{k-1} - t)(k-1)$ point-vertices can only be covered by the remaining $\frac{v}{k}-t-1$ cops, each of which is
adjacent to at most $k$ point-vertices.  Hence for the cops to surround the robber it is necessary that
$(\frac{v}{k}-t-1)k \geq (\frac{v-1}{k-1} - t)(k-1)$, which is clearly impossible.
Thus the robber is able to move away from vertex $x$, not be surrounded when the cops move, and then he can move to another point-vertex.
\end{proof}

Bonato and Burgess showed in~\cite{BonatoBurgess} that the cop number of the block-intersection graph of a BIBD$(v,k,\lambda)$ is at most $k$.
We establish the exact value for the surrounding cop number for BIBDs with index $\lambda=1$,
and in so doing we find that the surrounding cop number is substantially greater than the cop number.

\begin{theorem}
If $G$ is the block-intersection graph of a BIBD($v,k,1$), then $\sigma(G)=k(r-1)$.
\end{theorem}

\begin{proof}
First note that $\delta(G)=k(r-1)$ so $\sigma(G)\geq k(r-1)$. It remains to show that $\sigma(G)\leq k(r-1)$.  In the case that $G$ is a complete graph, then $\sigma(G)=\delta(G)=k(r-1)$, so we may assume that there are two non-intersecting blocks.  Hence we may
relabel the BIBD so that it contains blocks $\{1,2,\ldots, k\}$ and $\{k+1,k+2,\ldots,2k\}$.
Initially place the $k(r-1)$ cops on all blocks adjacent to $\{1,2,\ldots, k\}$. Without loss of generality, assume the robber starts on block $\{k+1,k+2,\ldots,2k\}$.

Let $A$ be the set of blocks adjacent to $\{1,2,\ldots, k\}$ in $G$, and let $B$ be the set of blocks adjacent to $\{k+1,k+2,\ldots,2k\}$ in $G$. As each pair of points appears together in a unique block, each $x\in A$ contains exactly one element in $\{1,2,\ldots, k\}$ and each $y\in B$ contains exactly one element in $\{k+1,k+2,\ldots, 2k\}$; thus $|A\cap B|=k^2$. We construct a bipartite graph $G'$ with partite sets $A'$ and $B'$ where $A'=\{x_A:x\in A\}$ is a copy of $A$ and $B'=\{y_B: y\in B\}$ is a disjoint copy of $B$. Let $x_Ay_B$ be an edge in $G'$ if $x_A\in A'$, $y_B\in B'$ and either $x$ is adjacent to $y$ in $G$ or $x=y$.

Let $H$ be the following subgraph of $G'$.
For each point $i\in\{1,2,\ldots,k\}$ include the edges $x_Ay_B \in E(G')$ such that  $i\in x$, $i+k \notin x$, and $i+k\in y$. Also include the edges $x_Ay_B$ where $\{i,i+k\}\subseteq x$ and $\{i+k,j\}\subseteq y$ for $j\in \{1,2,\ldots,k\}\setminus \{i\}$.
Clearly, $G'$ is a $(k-1)$-regular bipartite graph so there exists a $1$-factor $\mathcal{F}$ of $G'$.
For each edge $x_Ay_B$ in $\mathcal{F}$, the cops move from $x$ to $y$ in $G$, thereby surrounding the robber.
\end{proof}

\subsection{Generalised Petersen graphs}\label{Section-GeneralisedPetersen}

We begin this section with a definition.
Let $A=\{a_0, \ldots, a_{n-1}\}$ and $B=\{b_0, \ldots, b_{n-1}\}$.
For positive integers $n$ and $k$ such that $n > 2k$, the {\em generalised Petersen graph} $\GP(n,k)$ has vertex set $A \cup B$
and edges $\{a_i, a_{i+1}\}$, $\{a_i, b_i\}$ and $\{b_i, b_{i+k}\}$ for each $i \in \{0, \ldots, n-1\}$, with subscripts computed modulo $n$.
The edges of the form $\{a_i,b_i\}$ are called {\em spokes}.
The infinite Petersen graph $\GP(\infty,k)$ is defined similarly with $A=\{a_i:i\in\mathbb{Z}\}$ and $B=\{b_i:i\in\mathbb{Z}\}$.
Observe that the well known Petersen graph is $\GP(5,2)$.  In $\GP(n,k)$ or $\GP(\infty,k)$, we refer to a movement by a cop or robber as {\em leftwards} if the player moves from $a_i$ to $a_{i-1}$ or from $b_i$ to $b_{i-k}$, and {\em rightwards} if the player moves from $a_i$ to $a_{i+1}$ or from $b_i$ to $b_{i+k}$.

Tables~\ref{Table-cGenPet} and~\ref{Table-sGenPet} show the cop number and the surrounding cop number for several generalised Petersen graphs $\GP(n,k)$,
as computed by our implementation of the algorithm of Bonato and Chiniforooshan~\cite{BC2009} (for the cop number) and Algorithm~\ref{alg1} (for the surrounding cop number).

\begin{table}[]
\center
\begin{tabular}{|l|l|l|l|l|l|l|l|l|l|l|l|l|l|l|l|}
\hline
\multicolumn{2}{|l|}{\multirow{2}{*}{$c(\GP(n,k))$}} & \multicolumn{14}{l|}{\hspace{4.5cm}$k$}                                    \\ \cline{3-16}
\multicolumn{2}{|l|}{}                            & 1 & 2 & 3 & 4 & 5 & 6 & 7 & 8 & 9 & 10 & 11 & 12 & 13 & 14 \\ \hline
\multirow{28}{*}{$n$}              & 3              & 2 & &   &   &   &   &   &   &   &    &    &    &    &    \\ \cline{2-16}
                                 & 4              & 2 &   &   &   &   &   &   &   &   &    &    &    &    &    \\ \cline{2-16}
                                 & 5              & 2 & 3 &   &   &   &   &   &   &   &    &    &    &    &    \\ \cline{2-16}
                                 & 6              & 2 & 2 &   &   &   &   &   &   &   &    &    &    &    &    \\ \cline{2-16}
                                 & 7              & 2 & 3 & 3 &   &   &   &   &   &   &    &    &    &    &    \\ \cline{2-16}
                                 & 8              & 2 & 2 & 3 &   &   &   &   &   &   &    &    &    &    &    \\ \cline{2-16}
                                 & 9              & 2 & 3 & 2 & 3 &   &   &   &   &   &    &    &    &    &    \\ \cline{2-16}
                                 & 10             & 2 & 3 & 3 & 3 &   &   &   &   &   &    &    &    &    &    \\ \cline{2-16}
                                 & 11             & 2 & 3 & 3 & 3 & 3 &   &   &   &   &    &    &    &    &    \\ \cline{2-16}
                                 & 12             & 2 & 3 & 2 & 3 & 3 &   &   &   &   &    &    &    &    &    \\ \cline{2-16}
                                 & 13             & 2 & 3 & 3 & 3 & 3 & 3 &   &   &   &    &    &    &    &    \\ \cline{2-16}
                                 & 14             & 2 & 3 & 3 & 3 & 3 & 3 &   &   &   &    &    &    &    &    \\ \cline{2-16}
                                 & 15             & 2 & 3 & 3 & 3 & 3 & 3 & 3 &   &   &    &    &    &    &    \\ \cline{2-16}
                                 & 16             & 2 & 3 & 3 & 3 & 3 & 3 & 3 &   &   &    &    &    &    &    \\ \cline{2-16}
                                 & 17             & 2 & 3 & 3 & 3 & 3 & 3 & 3 & 3 &   &    &    &    &    &    \\ \cline{2-16}
                                 & 18             & 2 & 3 & 3 & 3 & 3 & 3 & 3 & 3 &   &    &    &    &    &    \\ \cline{2-16}
                                 & 19             & 2 & 3 & 3 & 3 & 3 & 3 & 3 & 3 & 3 &    &    &    &    &    \\ \cline{2-16}
                                 & 20             & 2 & 3 & 3 & 3 & 3 & 3 & 3 & 3 & 3 &    &    &    &    &    \\ \cline{2-16}
                                 & 21             & 2 & 3 & 3 & 3 & 3 & 3 & 3 & 3 & 3 & 3  &    &    &    &    \\ \cline{2-16}
                                 & 22             & 2 & 3 & 3 & 3 & 3 & 3 & 3 & 3 & 3 & 3  &    &    &    &    \\ \cline{2-16}
                                 & 23             & 2 & 3 & 3 & 3 & 3 & 3 & 3 & 3 & 3 & 3  & 3  &    &    &    \\ \cline{2-16}
                                 & 24             & 2 & 3 & 3 & 3 & 3 & 3 & 3 & 3 & 3 & 3  & 3  &    &    &    \\ \cline{2-16}
                                 & 25             & 2 & 3 & 3 & 3 & 3 & 3 & 4 & 3 & 3 & 3  & 3  & 3  &    &    \\ \cline{2-16}
                                 & 26             & 2 & 3 & 3 & 3 & 3 & 3 & 3 & 3 & 3 & 4  & 3  & 3  &    &    \\ \cline{2-16}
                                 & 27             & 2 & 3 & 3 & 3 & 3 & 4 & 3 & 3 & 3 & 3  & 3  & 3  & 3  &    \\ \cline{2-16}
                                 & 28             & 2 & 3 & 3 & 3 & 3 & 4 & 3 & 4 & 3 & 3  & 3  & 3  & 3  &    \\ \cline{2-16}
                                 & 29             & 2 & 3 & 3 & 3 & 3 & 3 & 3 & 4 & 3 & 3  & 4  & 4  & 3  & 3   \\ \cline{2-16}
                                 & 30             & 2 & 3 & 3 & 3 & 3 & 3 & 3 & 3 & 3 & 3  & 3  & 3  & 3  & 3  \\ \hline
\end{tabular}
\caption{Cop numbers of Generalised Petersen Graphs, $c(\GP(n,k))$}
\label{Table-cGenPet}
\end{table}

\begin{table}[]
\center
\begin{tabular}{|l|l|l|l|l|l|l|l|l|l|l|}
\hline
\multicolumn{2}{|l|}{\multirow{2}{*}{$\sigma(\GP(n,k))$}} & \multicolumn{9}{l|}{\hspace{2.5cm}$k$}            \\ \cline{3-11}
\multicolumn{2}{|l|}{}                            & 1 & 2 & 3 & 4 & 5 & 6 & 7 & 8 & 9 \\ \hline
\multirow{18}{*}{$n$}              & 3              & 3 &   &   &   &   &   &   &   &   \\ \cline{2-11}
                                 & 4              & 3 &   &   &   &   &   &   &   &   \\ \cline{2-11}
                                 & 5              & 3 & 3 &   &   &   &   &   &   &   \\ \cline{2-11}
                                 & 6              & 3 & 3 &   &   &   &   &   &   &   \\ \cline{2-11}
                                 & 7              & 3 & 4 & 4 &   &   &   &   &   &   \\ \cline{2-11}
                                 & 8              & 3 & 3 & 3 &   &   &   &   &   &   \\ \cline{2-11}
                                 & 9              & 3 & 4 & 3 & 4 &   &   &   &   &   \\ \cline{2-11}
                                 & 10             & 3 & 4 & 4 & 4 &   &   &   &   &   \\ \cline{2-11}
                                 & 11             & 3 & 4 & 4 & 4 & 4 &   &   &   &   \\ \cline{2-11}
                                 & 12             & 3 & 4 & 3 & 4 & 4 &   &   &   &   \\ \cline{2-11}
                                 & 13             & 3 & 4 & 4 & 4 & 4 & 4 &   &   &   \\ \cline{2-11}
                                 & 14             & 3 & 4 & 4 & 4 & 4 & 4 &   &   &   \\ \cline{2-11}
                                 & 15             & 3 & 4 & 4 & 4 & 4 & 4 & 4 &   &   \\ \cline{2-11}
                                 & 16             & 3 & 4 & 4 & 4 & 4 & 4 & 4 &   &   \\ \cline{2-11}
                                 & 17             & 3 & 4 & 4 & 4 & 4 & 4 & 4 & 4 &   \\ \cline{2-11}
                                 & 18             & 3 & 4 & 4 & 4 & 4 & 4 & 4 & 4 &   \\ \cline{2-11}
                                 & 19             & 3 & 4 & 4 & 4 & 4 & 4 & 4 & 4 & 4 \\ \cline{2-11}
                                 & 20             & 3 & 4 & 4 & 4 & 4 & 4 & 4 & 4 & 4 \\ \hline
\end{tabular}
\caption{Surrounding Cop numbers of Generalised Petersen Graphs, $\sigma(\GP(n,k))$}
\label{Table-sGenPet}
\end{table}

The cop number of generalised Petersen graphs was previously considered by
Ball, Bell, Guzman, Hanson-Colvin and Schonsheck in~\cite{Ball2017}.
They proved that every generalised Petersen graph has cop number at most 4, and by performing a computational assessment they demonstrated that
a cop number of 4 is required for several cases.
With the exception of $\GP(25,7)$, our computational results agree with theirs.
Our implementation has stated that $c(\GP(25,7)) = 4$ (see Table~\ref{Table-cGenPet}),
whereas Ball et al.\ reported that $\GP(26,10)$ is the smallest generalised Petersen graph having 4 as its cop number.

In Theorem~\ref{Lemma-GeneralisedPetersen} we show that every generalised Petersen graph has surrounding cop number at most 4.
Empirical results shown in Table~\ref{Table-sGenPet} suggest that there may only be a handful of cases for which
$\sigma(\GP(n,k)) = 3$ and $k \geq 2$.

The proof of Theorem~\ref{Lemma-GeneralisedPetersen} is similar to that of~\cite{Ball2017} for ordinary Cops and Robbers.  The key is to consider a ``lifting'' of the game on GP$(n,k)$ to a game on the infinite graph GP$(\infty,k)$.  Before proceeding to the proof, we require some further terminology adapted from \cite{Ball2017}.

Consider the homomorphism $\pi: \GP(\infty,k) \rightarrow \GP(n,k)$ which maps a vertex $x_i$ (where $x \in \{a,b\}$) to the vertex $x_j$ of GP$(n,k)$ satisfying $j \equiv i$ (mod $n$). If a player moves from vertex $y$ to $z$ in GP$(n,k)$, the {\em lifted move} in GP$(\infty,k)$ for a player on a vertex $v$ of $\pi^{-1}(y)$ is to move to the unique vertex of $\pi^{-1}(z)$ adjacent to $v$, or to pass if $y=z$.  Suppose in GP$(\infty,k)$, two players, say $P_1$ and $P_2$, are on vertices $u$ and $v$, respectively, of $\pi^{-1}(y)$.  We describe the next move of $P_1$ and $P_2$ as {\em consistent} if $P_1$ and $P_2$ follow the lifted move of the same move in GP$(n,k)$.  Given a game of Surrounding Cops and Robbers in GP$(n,k)$, the {\em lifted game} in GP$(\infty,k)$ is as follows.  For each cop $\cal{C}$ in GP$(n,k)$, we associate a set $\cal{S}$ of cops in GP$(\infty,k)$; if $\cal{C}$ is initially placed on vertex $x$, then a cop of $\cal{S}$ is initially placed on each vertex of $\pi^{-1}(x)$.  Thereafter, in each round the cops of $\cal{S}$ play the lifted moves of those of $\cal{C}$; we refer to the cops of $\cal{S}$ as a {\em squad}.  Suppose the robber $R$ begins on vertex $y$ of GP$(n,k)$.  We place a robber on a vertex of $\pi^{-1}(y)$, who will play the lifted moves of $R$.  Conversely, if a game in GP$(\infty,k)$ is the lift of a game in GP$(n,k)$, we refer to the game in GP$(n,k)$ as the {\em projected game}.

In practice, while playing in GP$(\infty,k)$, we single out one cop of a given squad, the {\em lead cop}, and determine the consistent movements of the entire squad based on the position of the lead cop relative to the robber.  As the game progresses, the choice of lead cop may change; however, regardless of this choice, the movements of the squad correspond to the movements of a single cop in the projected game.

\begin{lemma}\label{Lemma-GPliftedgame}
In the game of Surrounding Cops and Robbers, two squads of cops playing in $\GP(\infty,k)$ can prevent the robber from moving left after a finite number of turns.
\end{lemma}
\begin{proof}
The proof is similar to that of Corollary 2.2 of~\cite{Ball2017}, with the exception that the robber cannot necessarily be forced to increase his index, but rather is prevented from decreasing his index beyond a finite number of moves.  We provide only an outline of the details. The difference in the cops' strategy is that the cops now choose to pass rather than landing on the robber's position.

Using the modified strategy,
we outline how, after finitely many turns, the robber is prevented from moving left.
The cops start on $A$ and work in tandem so that
the first lead cop, ${\cal C}_1$, matches parity and congruence modulo $k$ with the robber.
In subsequent rounds, the cop ${\cal C}_1$ then moves (if necessary) to remain in the same set ($A$ or $B$) as the robber and maintain congruence modulo $k$, and
each time the robber moves within $B$ or passes while on a vertex of $B$, ${\cal C}_1$ moves closer to the robber in $B$.  Thus, the robber can only move or pass finitely many times within $B$ before ${\cal C}_1$ and the robber occupy vertices of the form $x_i$ and $x_{i+k}$ ($x\in \{a,b\})$; if this happens their indices will differ by $k$ for the remainder of the game, thus preventing the robber from making any leftward moves in $B$.  However, each time the robber moves leftwards within $A$, the second lead cop (${\cal C}_2$) moves closer to the robber in $A$, so the robber can make only finitely many leftward moves in $A$.  (Note that if the robber moves along a spoke, ${\cal C}_2$ will not move to a vertex with larger index than the robber, so movement along spokes will not increase the robber's number of leftward moves in $A$.)
\end{proof}

In~\cite{Ball2017} it was proved that the cop number of every generalised Petersen graph is at most 4.  This constant upper bound also applies to the surrounding cop number.

\begin{theorem}\label{Lemma-GeneralisedPetersen}
For all integers $n\geq 5$ and $k\geq 2$, $\sigma(\GP(n,k))\leq 4$.
\end{theorem}

\begin{proof}
The proof due to Ball et al.~\cite{Ball2017} can be adapted for $\sigma(G)$. This proof relies on playing the lifted version of the game on the graph $\hat{G}=\GP(\infty,k)$. The main change is that we use Lemma~\ref{Lemma-GPliftedgame} in place of~\cite[Corollary 2.2]{Ball2017}.

By Lemma~\ref{Lemma-GPliftedgame}, one pair of lifted cops can prevent $\hat{R}$ (the lifted robber) from moving to the left. A second pair of lifted cops moves closer to $\hat{R}$ or prevents $\hat{R}$ from moving to the right.  At this point, the robber's movements in the lifted game, and hence in GP$(n,k)$, are confined to the spokes, and the cops are in position to surround him by the next turn.
\end{proof}

As first noted in~\cite{Ball2017}, there are many graphs $\GP(n,k)$ where $c(\GP(n,k)) = 3$ and seemingly many others for which $c(\GP(n,k)) = 4$.
Empirical results (see Table~\ref{Table-sGenPet})
suggest that when $k \geq 2$ there may only be a finite number of instances for which the surrounding cop number is not 4,
namely $\GP(5,2)$, $\GP(6,2)$, $\GP(8,2)$, $\GP(8,3)$, $\GP(9,3)$ and $\GP(12,3)$.
In the case $k=1$, $\GP(n,1)$ is a prism and $\sigma(\GP(n,1)) = 3$.

\begin{theorem}
For all $n \geq 3$, $\sigma(\GP(n,1)) = 3$.
\end{theorem}

\begin{proof}
Since $\GP(n,1)$ is 3-regular, clearly at least three cops are needed to surround the robber.  We give a strategy in which three cops can prevail.  We let $C_A$ (resp.\ $C_B$) denote the cycle induced by $A$ (resp.\ $B$).

Initially, the cops are placed with two of them, say $\cC_1$ and $\cC_2$, on $a_0$.  The third, $\cC_3$, is placed on $b_0$.  In the first stage of the game, cops $\cC_1$ and $\cC_2$ play on the cycle induced by $A$, aiming to surround the robber on the graph obtained by identifying vertex $a_i$ with $b_i$ for each $i$.  To simplify the discussion below, define the {\emph{shadow}} of the robber to be $a_i$ if the robber is on vertex $b_i$ and $b_i$ if the robber is on $a_i$.  In this part of the game, $\cC_3$ passes in each round, thus preventing the robber from moving cyclically along $C_B$ (this ensures that if the robber is playing on $C_B$, cops $\cC_1$ and $\cC_2$ can surround his shadow). At the end of this stage, we have one of two cases.

\vspace*{0.5ex}
\noindent
{\bf Case 1.}  $\cC_1$ and $\cC_2$ surround the robber on $C_A$, so that $\cC_1$ and $\cC_2$ are in positions $a_{i-1}$ and $a_{i+1}$, with the robber on vertex $a_i$.

In this case, $\cC_1$ and $\cC_2$ confine the robber to moving along the spoke $a_i b_i$.  If he passes, so do they; otherwise, his only allowable move is along this spoke, and $\cC_1$ and $\cC_2$ follow along their corresponding spokes.  At this stage of the game, $\cC_3$ moves until she is at vertex $b_i$.  If the robber is on $a_i$ when $\cC_3$ arrives at $b_i$, the robber is immediately surrounded.  Otherwise, $\cC_3$ lands on the robber, who is then forced to move to $a_i$.  In the next move, cops $\cC_1$ and $\cC_2$ (who had been on $b_{i-1}$ and $b_{i+1}$) move to $a_{i-1}$ and $a_{i+1}$, while $\cC_3$ passes; again, the robber is surrounded.

\vspace*{0.5ex}
\noindent
{\bf Case 2.}  $\cC_1$ and $\cC_2$ surround the robber's shadow, that is, they occupy positions $a_{i-1}$ and $a_{i+1}$ while the robber is on vertex $b_i$.

If the robber moves to $a_i$, we proceed as in Case 1.  Otherwise, after the next move of the cops or the robber, we can ensure that one of the cops occupies the other end-vertex of the spoke incident with the robber.  In the rest of the game, this cop stays on $C_A$, mirroring the robber's moves on $C_B$, so that after each of this cop's turns she is still on the other end of the robber's spoke, and the robber is henceforth prevented from moving to a vertex of $A$. For the robber and the other two cops, the game is reduced to playing on the cycle $C_B$, where two cops can surround the robber.
\end{proof}

\section{Concluding Remarks}

Although the game of Surrounding Cops and Robbers is new, it has already attracted interest.  A recent preprint~\cite{BradshawHosseini} considers the surrounding cop number of graphs of bounded genus, in particular giving small constant upper bounds for the surrounding cop number of planar, outerplanar and toroidal graphs. Below we present several unresolved questions to provide directions for future study.

\begin{enumerate}
\item The empirical results that we have for generalised Petersen graphs suggest that $\sigma(\GP(n,k))=4$ whenever $n > 12$ and $k > 1$.  Is this indeed true?
\item In~\cite{CKM2020} it is proved that $c(G) \leq \xi(G)$.
Is the containability number $\xi(G)$ also an upper bound for $\sigma(G)$?
\item Theorem~\ref{Theorem-Girth} establishes $\sigma(G) > \delta$
whenever $G$ has minimum degree $\delta \geq 3$ and girth $g \geq 7$.
Does a generalisation of this theorem, similar to the result of Frankl~\cite{Frankl1987}, hold?
\item Characterise graphs $G$ such that $\sigma(G)=k$.
For $k=1$ we found that $G$ must be $K_1$ or a star $K_{1,m}$ for some $m \geq 1$.
For $k=2$ we do not have a characterisation, but have noted that the family includes such graphs as cycles, trees (other than $K_1$ and stars),
cycles with a chord, and more.
\end{enumerate}

It is not true that $G$ being a connected subgraph of a connected graph $H$  implies that $\sigma(G) \leq \sigma(H)$.
For instance, form the graph $G$ from an 8-cycle $(1,2,3,\ldots,8)$ by adding the edges $\{2,8\}$ and $\{4,6\}$ and let $H = G+e$ where $e$ is the edge $\{2,6\}$. It is easy to confirm that $\sigma(G) = 3$ and $\sigma(H) = 2$.

\begin{enumerate}\setcounter{enumi}{4}
\item Under what conditions on $G$ and $e$ is $\sigma(G-e)\leq \sigma(G)$?
\item Let $G/e$ be the graph obtained by contracting the edge $e$.  Under what conditions on $G$ and $e$ is $\sigma(G/e)\leq \sigma(G)$?
\end{enumerate}

We conclude with a brief discussion regarding Meyniel's conjecture, which asserts that
$c(G) \in \bigO(\sqrt{n})$, where $n$ denotes the order of the graph $G$. This conjecture was first published in \cite{Frankl1987}.
For a nice survey, see \cite{bb}.
When considering graphs such as complete graphs, for which $\sigma(K_n)=n-1$,
it is clear that there is no analogy of this conjecture that would apply to the surrounding cop number. Nevertheless, there may be a way to exploit the surrounding cop number
as a means of proving Meyniel's conjecture for the cop number.
Recall that $c(G) \leq \sigma(G)$, and in particular note that any graph $G$ for which $\sigma(G) \in \bigO(\sqrt{n})$ must satisfy Meyniel's conjecture.
Hence any counterexample to the conjecture must be such that $\sigma(G) \notin \bigO(\sqrt{n})$.
Aside from incidence graphs of projective planes,
all of the graphs for which $\sigma(G)$ exceeds $\sqrt{n}$ that we have encountered have been graphs
for which the cop number $c(G)$ has not only been below $\sqrt{n}$, but substantially below.
For instance, when $G$ is an $n$-vertex line graph of a complete graph, in which case $\sigma(G)$ is approximately $2\sqrt{2n}$ (see Theorem~\ref{Thm-LineGraphKn}),
the cop number is at most 2.
\begin{enumerate}\setcounter{enumi}{6}
\item Do graphs with high surrounding cop number 
inherently possess some property which in turn implies that the cop number is low?
\end{enumerate}

\section{Acknowledgements}
Authors Burgess, Clarke, Danziger, Finbow and Pike acknowledge research grant support from NSERC discovery grants RGPIN-2019-04328, RGPIN-2015-06258, RGPIN-2016-04178, RGPIN-2014-06571 and RGPIN-2016-04456, respectively.
Cameron acknowledges support from an AARMS postdoctoral fellowship.
Jones acknowledges support from an NSERC Undergraduate Student Research Award.

\end{document}